\documentstyle[12pt]{amsart}   
  \title[$\mbox{\em Fixed points of elliptic reversible transformations}$] {$\mbox{\Huge\bf Fixed points of elliptic reversible}$\\ 
$\mbox{\Huge\bf transformations with 
integrals}$}
\thanks{{\it Current address:} Institute for Advanced Study, School of Mathematics, Olden Lane, Princeton, NJ 08540, USA}
 \author[$\mbox{\em X$.$ Gong}$]{ Xianghong  Gong}
 \dedicatory{
Department of Mathematics$,$  University of Chicago$,$
Chicago$,$ IL $60637,$ USA$\dag$}

 \newtheorem{thm}{Theorem}[section]
 \newtheorem{cor}[thm]{Corollary}
 \newtheorem{prop}[thm]{Proposition}

 \newtheorem{lemma}[thm]{Lemma}
 
\begin{document}
 \maketitle
 \vspace{2ex} \noindent{\it Abstract}.
We show that for a certain
family of integrable reversible transformations, the curves of periodic
points of a general transformation  cross the level curves of 
its integrals. This leads to the divergence of the normal form for a general 
reversible transformation with integrals. We also study the integrable
holomorphic reversible transformations coming from real analytic surfaces
in ${\Bbb C}^2$
with non-degenerate complex tangents. We show the existence of real
analytic surfaces with hyperbolic complex tangents, which are contained in
a real
hyperplane, but cannot be
transformed into the Moser-Webster normal form through any holomorphic
transformation.

\setcounter{section}{1}
\vspace{3ex}\noindent
{1.\hspace{.5em}\it Introduction and results}\\
\noindent
A transformation $\varphi\colon {\Bbb R}^2\to{\Bbb R}^2$ is {\it reversible}
if it is  conjugate to its inverse by  an involution, or equivalently if $\varphi$ is the composition of a pair of involutions. 
 Reversible transformations  appeared often in
Birkhoff's work on dynamical systems due mainly  to the existence of periodic
points of  a general elliptic reversible transformation
(see~\cite{birkhoff1},~\cite{birkhoff3}). The  periodic
points of reversible mappings and   systems were further 
discussed by R.\ ~L.\ ~Devaney~\cite{devaney}. It
has been observed that reversible mappings and reversible systems  have many properties similar to
area-preserving transformations or Hamiltonian systems~\cite{moser5},~\cite{arnold2}.
For instance, using the curve intersection property, J.\ ~K.\ ~Moser~\cite{moser4} showed that  area-preserving
 mappings  have invariant curves surrounding
 elliptic fixed points. For reversible transformations, the existence of 
invariant curves was first announced by
V.\ ~I.\ ~Arnol'd in~\cite{arnold2}. A complete proof was given by  M.\ ~B.\ ~Sevryuk in~\cite{sevryuk} where examples of  reversible transformations which
do not
satisfy the curve intersection property were also constructed.

The purpose of this paper is to  bring up a different aspect of reversible
transformations concerning  the convergence of normal forms.  In~\cite{birkhoff2},
G.\ ~D.\ ~Birkhoff proved that an area preserving transformation $\psi\colon {\Bbb R}^2\to{\Bbb R}^2 $ can always be transformed into its normal form through convergent transformations  whenever $\psi$  has an integral, i.\ ~e.\  a non-constant function invariant under $\psi$.  However, we shall prove that an analytic reversible transformation may not be transformed into its normal form by any
 convergent transformation even if it has  integrals.

 Let $z=x+iy$ be the complex coordinate for $R^2$, and  denote by $\tau$ the complex conjugation $z\to\overline z$. To formulate our result,  we fix a positive
integer $s$ and let 
$$\sum=\{ a=(a_{i,j});i+j>2s, a_{j,i}=\overline a_{i,j},|a_{i,j}|<1\}.$$
We shall also regard $a_{i,j}$ as a variable and denote by $a$ the infinitely many variables $a_{i,j}\ (i+j>2s)$.
Thus,  $h(a)$ will stand for a power series in all variables $a_{i,j} $. For each $a\in\sum$, we put
$$\widetilde a(z,\overline z)=
 \sum_{i+j>2s}a_{i,j}z^i\overline z^j.
 $$
We define
 a mapping $\varphi_a\colon{\Bbb R}^2\to{\Bbb R}^2$ by
 \begin{gather}\label{eq:1.3}
 z\to \phi_aT\tau\phi_a^{-1}\tau(z), \\
\intertext{where $\phi_a(z)=ze^{i\widetilde a(z,\overline z)}$, and $T$ is the twist mapping }
z\to ze^{i(\alpha+(z\overline z)^s)},\quad \alpha/(2\pi)\in{\Bbb R}\setminus {\Bbb Q}.\label{eq:1.2}
\end{gather}
We shall see that $T$ is the formal normal form of  $\varphi_a$. Notice that
$\varphi_a$ is reversible with respect to $\tau$, i.\ ~e.\
$\varphi_a^{-1}=\tau\varphi_a\tau^{-1}$.
 Moreover,  the function
$\kappa(z,\overline z)=z\overline z$ is an integral for all $\varphi_a$.
 
 We shall prove the following.
\begin{thm}\label{cor:1.2}
{\it For each $\alpha/(2\pi)\notin {\Bbb Q},$ there exists  a sequence of functionally independent   power series
 $H_k$ in $a,$ of which each $H_k$ converges for $a\in\sum,$
 such that for a fixed $a\in\sum,$ $\varphi_a$  cannot be
 transformed into $T$
 through any convergent transformation of ${\Bbb R}^2,$ if there are infinitely many $H_k(a)\neq
0.$
}
 \end{thm} 
 Theorem~\ref{cor:1.2} is motivated by a result of C.\ ~L.\ ~Siegel about Hamiltonian systems~\cite{siegel2} as well as a result of 
 of H.\ ~R\"{u}ssmann about area-preserving
transformations~\cite{russmann3}. We shall prove Theorem~\ref{cor:1.2} by  investigating the periodic points of $\varphi_a$. Our observation is that as a  perturbation of the twist mapping,  $\varphi_a$ has  the Birkhoff curves surrounding the origin.
In fact, these curves are the curves of periodic points since $\varphi_a$ preserves all the circles centered at the origin. However, we shall prove that, as a rule,  the  curves of periodic points of $\varphi_a$ are not contained in any circle centered at the origin. Thus, we have the following result.
 \begin{thm}\label{thm:1.1}
{\it
For each $\alpha/(2\pi)\notin {\Bbb Q},$ there exists  a sequence of functionally independent convergent power series 
 $H_k(a), a\in\sum,$ which satisfies the following property$:$ If $a\in\sum$ is fixed and  $H_k(a)\neq 0$ for some $k,$
then there is a closed interval $I_k$ of positive length 
such that for each $r\in I_k,$  $\varphi_a$ has a periodic point on $|z|=r.$ Furthermore$,$ $I_k$ is contained in any given interval $(0,\epsilon)$ when $k$ is large$.$
}
 \end{thm}
Notice that there are only countable many circles centered at the orign on which the
restriction of $T$ has  a periodic point. Thus,
Theorem~\ref{cor:1.2} follows  Theorem~\ref{thm:1.1}.

 To state our next result, we consider a real analytic   surface in ${\Bbb C}^2$ with a non-degenerate complex tangent at the origin given by
 \begin{equation}\label{eq:1.4}
 M\colon z_2=z_1\overline z_1+\gamma z_1^2+\gamma \overline z_1^2+q(z_1,\overline z_1),\quad \gamma\in{\Bbb R},
 \end{equation}
 in which $q$ is a convergent power series starting with terms of order three.
When $0<\gamma<1/2$, or $1/2<\gamma<\infty$ with a countable 
exceptional values,  J.\ ~K.\ ~Moser and S.\ ~M.\ ~Webster~\cite{moserwebster}
showed that through formal transformations, $M$ can be transformed into a surface defined by
\begin{equation}\label{eq:norm}
  x_2=z_1\overline z_1+(\gamma +\epsilon x_2^s)(z_1^2+\overline z_1^2),\quad y_2=0,
\end{equation}
 where either $\epsilon=0$, or $\epsilon=\pm 1$ with $s$ a positive integer.
The normal form (\ref{eq:norm}) was constructed through a pair of
holomorphic involutions   which
characterizes the real analytic surfaces in ${\Bbb C}^2$ with a non-degenerate complex tangent. The existence of an integral  for such a pair of involutions corresponds to the property that the surface can be transformed into a real hyperplane in ${\Bbb C}^2$.  Examples of real analytic surfaces of hyperbolic complex tangents which cannot be transformed into any real hyperplane  were constructed in~\cite{moserwebster} and also in~\cite{bedford} by E.\ ~ Bedford.
Analogous to Theorem~\ref{cor:1.2}, We have the following.
 \begin{thm}\label{thm:1.3}{\it
 For each non-exceptional\/ $\gamma\in (1/2,\infty),$ there exists a
real analytic surface  $(\ref{eq:1.4})$  which can be transformed  biholomorphically into a real hyperplane in ${\Bbb C}^2,$ but  not into the normal form $(\ref{eq:norm}).$}
  \end{thm}
 We mention that for a set of $\gamma$ with Lebesgue measure $0$, the above
result was proved in~\cite{gong}.

The paper is organized as follows. In section $2$, we discuss the normal
form for reversible transformations. Theorem~\ref{thm:1.1} is proved in
section $4$ following the preparation in section $3$ where   the existence of
periodic points of holomorphic transformations with integrals is
discussed. Section $6$ is devoted to the study of holomorphic reversible
transformations in ${\Bbb C}^2$ which come from real surfaces with complex tangents. The proof of
Theorem~\ref{thm:1.3} is given in section $6$.

\setcounter{section}{2}
\vspace{3ex}\noindent
{2.\hspace{.5em}\it Normal form for reversible transformations}\\
\noindent
 \setcounter{thm}{0}\setcounter{equation}{0}
 Let $\varphi,\tau$ be two real analytic 
transformations of ${\Bbb R}^2$ defined near the origin. 
Assume that $\varphi(0)=\tau(0)=0$, and $\tau^2=id$. Suppose that
$\varphi$ is reversible with respect to $\tau$.  
With a suitable assumption on the eigenvalues of  $d\varphi(0)$, 
we shall see that after a change of analytic coordinates, $\tau$ is precisely the complex conjugation $z\to\overline z$.  Then we shall derive a normal form for $\varphi$ under the formal transformations which commute with $\tau$. The normal formal is essentially due to Moser and  Webster~\cite{moserwebster}.
 
  For convenience,  we shall complexify real analytic or formal transformations of ${\Bbb R}^2$. Let  $\Phi$ be  a real analytic transformation  of ${\Bbb R}^2$. Then  $\Phi$ is a power series in $z$ and $\overline z$. The complexification 
of $\Phi$ defined by
 $$
 (\xi,\eta)\to\bigl(\Phi(\xi,\eta),\overline\Phi( \eta, \xi)\bigr).
 $$
 is a holomorphic
transformation of ${\Bbb C}^2$, which will  still be denoted by  $\Phi$. Clearly,
 $\Phi$ satisfies the { \it reality condition }
 \begin{equation}\label{eq:2.1}
 \rho\Phi=\Phi\rho, 
 \end{equation}
 in which $\rho$ is the complexification of the complex conjugation, i.\ ~e.
 \begin{equation}\label{eq:2.2}  
 \rho:(\xi,\eta)\to(\overline \eta, \overline \xi).
 \end{equation}
 The complexfication of a fomal transformation is defined in a similar way. From now on, we shall identify the real analytic transformations of ${\Bbb R}^2$ with the  holomorphic transformations of ${\Bbb C}^2$ satisfying the reality condition (\ref{eq:2.1}).

 To discuss the normal form of $\varphi$, we first assume that the eigenvalues 
of $d\varphi(0)$ are not  roots of unity. Let $\lambda$ be an eigenvalue of  $d\varphi(0)$ with an eigenvector $e_1$.   We have
 \begin{equation}\label{eq:2.3}  
 d\varphi(0)\bigl(d\tau(0)e_1\bigr)=d\tau(0)(d\varphi^{-1}(0)e_1)=
 \lambda^{-1}d\tau(0)e_1.
 \end{equation}
 This implies that $\lambda^{-1}$ is also an eigenvalue of $d\varphi(0)$. Since 
 $\lambda^{-1}\neq\lambda$, we see that $d\varphi(0)$ is diagonalizable with distinct eigenvalues
 $\lambda$ and $\lambda^{-1}$. 

We now further assume that $\varphi$ is {\it elliptic\/}, i.\ ~e.\ ~$|\lambda|=1$ and $\lambda\neq\pm 1$. Put $e_2=\rho(e_1)$. Then the reality condition
 $\rho\varphi=\varphi\rho$ gives
 $$
 d\varphi(0)e_2=\rho d\varphi(0)\rho e_2=\overline \lambda e_2.
 $$
Hence, $\overline\lambda=\lambda^{-1}$. 
Under  new  coordinates $(\xi,\eta)$ for $\xi e_1+\eta e_2\in
{\Bbb C}^2$, 
the anti-holomorphic involution $\rho$ still takes the form (\ref{eq:2.2}),
while the reversible transformation is 
 \begin{equation}\label{eq:2.4}
 \varphi:\begin{split}&
 \xi^\prime=\lambda\xi+p(\xi,\eta),\\ &
 \eta^\prime=\overline \lambda\eta+q(\xi,\eta),
 \end{split}
 \end{equation}
 in which $p $ and $q $ are convergent power series starting with the second order
terms. Notice that $\lambda$ is replaced by $\overline \lambda$, if one
applies a change of coordinates by $(\xi,\eta)\to (\eta,\xi)$. Therefore,
one may assume that 
\begin{equation}\label{lambda}
\mbox{Im}\, \lambda>0.
\end{equation} 

We now want to normalize $\tau$. From  (\ref{eq:2.3}), we see that $d\tau(0)e_1$ is an eigenvector. Hence,
 $
 d\tau(0)e_1=\lambda_0 e_2.
 $
On the other hand, $\tau$ satisfies the reality condition. Then
 $$
 d\tau(0)e_1=\rho d\tau(0)\rho e_1=\rho d\tau(0)e_2=\overline \lambda_0^{-1}e_2.
 $$
 This implies that $|\lambda_0|=1$. Thus, we can write
 $$
 \tau(\xi,\eta)=(\lambda_0\eta,\overline \lambda_0\xi)+O(2).
 $$
 Consider a change of coordinates defined by
 \begin{equation*}
\begin{split}&
 \xi^\prime=\lambda_0^{-1/2}\bigl(\xi+\lambda_0\eta\circ\tau(\xi,\eta)\bigr)/2,\\ &
 \eta^\prime=\lambda_0^{1/2}\bigl(\eta+\overline \lambda_0\xi\circ\tau(\xi,\eta)\bigr)/2.\end{split}
 \end{equation*}
 From the reality condition $\rho\tau=\tau\rho$, one can see that under new coordinates $(\xi^\prime,\eta^\prime)$, $\rho$ is still of the form (\ref{eq:2.2}). Now $ \tau $ is the linear
 involution 
 \begin{equation}\label{eq:2.5}
 (\xi,\eta) \to(\eta,\xi).
 \end{equation}
 \begin{thm}\label{thm:2.1}{\it
 Let $\rho$ and $\tau $  be given by $(\ref{eq:2.2})$ and
 $(\ref{eq:2.5})$ respectively$.$ Suppose that $\varphi$ defined by $(\ref{eq:2.4})$
 is  reversible with respect to $\tau, $ and  it satisfies
 the reality condition $\rho\varphi=\varphi\rho.$  If  $\lambda$ is not a root of unity$,$ then there exists a formal transformation
 $\Phi$ such that $\rho\Phi=\Phi\rho, \tau \Phi=\Phi\tau $ and
 \begin{equation}\label{eq:normalform}
\varphi^*= \Phi\varphi\Phi^{-1}\colon\begin{split}&
 \xi^\prime=\lambda\xi e^{i\epsilon(\xi\eta)^s},\\ &
 \eta^\prime=\overline\lambda\eta e^{-i\epsilon(\xi\eta)^s},
 \end{split}
 \end{equation}
 in which $\epsilon=\pm 1$ with  $s$  a positive integer$,$ or
 $\epsilon=0$ with $s=\infty.$ Furthermore$,$ $\{\lambda, \epsilon,s\}$ is the
 full set of invariants of $\varphi$ under real formal transformations$.$
}
 \end{thm}

The  normal form (\ref{eq:normalform}) was derived
in~\cite{moserwebster} except a different reality condition is involved here.
For completeness, we shall verify the reality condition here. A transformation 
\begin{equation}\label{eq:2.16}
\xi\to\xi+u(\xi,\eta),\quad \eta\to\eta+v(\xi,\eta), \quad u,v=O(2)
\end{equation}
is said to be {\it normalized\/} if  $u_{i+1,i}=v_{i,i+1}=0$ for all $i$.
We need the following.
\begin{lemma}[\cite{moserwebster}]\label{lemma:2.2}{\it 
Let $\tau_j\  (j=1,2)$ be a pair of involutions 
\begin{equation}\label{eq:2.15}
 \xi\to \lambda_j\eta+f_j(\xi,\eta),\quad \eta\to \lambda_j^{-1}\xi +g_j(\xi,\eta).
 \end{equation}
Put $\varphi=\tau_1\tau_2.$ If $\lambda_1\lambda_2^{-1}$  is not a root of unity$,$ then there exists a unique normalized formal transformation $\Phi$    such that $\Phi\varphi\Phi^{-1}$ is of the form
\begin{equation}\label{eq:2.17}
\widetilde \varphi\colon \xi\to M(\xi\eta)\xi,\quad \eta\to M^{-1}(\xi\eta)\eta
\end{equation}
with $M(t)=\lambda_1\lambda_2^{-1}+O(1).$ Furthermore$,$ $\tau_j$ is transformed into 
\begin{equation}\label{eq:2.18}
\widetilde\tau_j\colon \xi\to\Lambda_j(\xi\eta)\eta,\quad \eta\to\Lambda_j^{-1}(\xi\eta)\xi
\end{equation}
with $\Lambda=\lambda_j+O(1).$
}
\end{lemma}

We now turn to the proof of Theorem~\ref{thm:2.1}. We put $\tau_1=\tau$ and
$\tau_2=\tau\varphi$. Then for the pair of involutions $\{\tau_1,\tau_2\}$,
it follows from Lemma~\ref{lemma:2.2} that there is a
unique normalized transformation $\Phi_0$   which transforms $\varphi$ into  (\ref{eq:2.17}).  Since $\rho\varphi=\varphi\rho$, we get
$$
(\rho\Phi_0\rho)\varphi(\rho\Phi_0\rho)^{-1}(\xi,\eta)=\bigl(\overline
M^{-1}(\xi\eta) \xi,\overline M(\xi\eta)\eta\bigr).
$$
Clearly, $\rho\varphi\rho$ is still in the form (\ref{eq:2.16}) except
$\overline u,\overline v$ are replaced by  $u$ and $v$, respectively. By the uniqueness of 
the transformation  $\Phi_0$, we have
\begin{equation}\label{eq:5.19}
\rho\Phi_0\rho=\Phi_0,\quad \overline M^{-1}(\xi\eta)=M(\xi\eta).
\end{equation}
Hence, $\rho\widetilde\tau_j=\widetilde\tau_j\rho$. This implies that
$\overline\Lambda_j^{-1}(t)= \Lambda_j(t)$. Put $a(t)=\Lambda_1^{1/2}(t)$ 
and 
$$
\Phi_1(\xi,\eta)=\bigl(a(\xi\eta)\xi,a^{-1}(\xi\eta)\eta\bigr).
$$
Then $\Phi_1^{-1}(\xi,\eta)=(a^{-1}(\xi\eta)\xi,a(\xi\eta)\eta)$.
 It is easy to see that
$$
\Phi_1\widetilde\tau_1\Phi_1^{-1}=\tau,\quad \Phi_1\rho=\rho\Phi_1,\quad
\Phi_1\widetilde\varphi\Phi_1^{-1}=\widetilde\varphi.
$$
We now write
$$
M(t)=e^{i\Gamma(t)}.
$$
Then (\ref{eq:5.19}) gives $\Gamma(\xi\eta)=\overline\Gamma(\xi\eta)$.
If $\Gamma(t)$ is not constant, we can find a real
power series $r(t)$ such that
$$
\Gamma(t)=\epsilon t^sr^{2s}(t),
$$
where $\epsilon$, $s$ are the sign and the order of the first non-vanishing
coefficient of $\Gamma(t)$, respectively.
Let $\Phi_2(\xi,\eta)=(\xi r(\xi\eta),\eta r(\xi\eta))$.
 Then one can verify that
$$
\Phi_2\tau=\tau\Phi_2,\quad \Phi_2\rho=\rho\Phi_2.
$$
Now $\Phi_2$ transforms $\widetilde\varphi$ into (\ref{eq:normalform}). 
 Therefore, $\Phi=\Phi_2\Phi_1\Phi_0$ preserves $\tau$ and $\rho$, and it
transforms $\varphi$ into (\ref{eq:normalform}). 

To show that $\epsilon$ and $s$ are invariants, we
assume that there is a formal transformation $\Psi$ such that
\begin{equation}\label{eq:5.20}
\Psi\varphi^*=\psi\Psi,
\end{equation}
where $\psi$ is a transformation in the form (\ref{eq:normalform}), of which
$\epsilon$ and $s $ are replaced by $\epsilon^\prime$ and $s^\prime$. We
may assume that $s\leq s^\prime$. From (\ref{lambda}), we first notice that
$\lambda$ is an invariant and 
$d\Psi(0)=(a_0\xi,b_0\eta)$. 
 Notice that $\lambda$ is not a root of unity. By comparing  both sides of  (\ref{eq:5.20}) up to terms of order
 $2s$,
 one gets
$$
\Psi(\xi,\eta)=\bigl(\xi a(\xi\eta),\eta b(\xi\eta)\bigr)+O(2s+2).
$$
  Now the reality condition
(\ref{eq:2.1}) implies that $b=\overline a$.
By comparing terms in (\ref{eq:5.20}) of  order $2s+1$,
one can get $s=s^\prime$ and $\epsilon=\epsilon^\prime$.
This proves Theorem~\ref{thm:2.1}.

Let
$T$ be the  twist mapping   (\ref{eq:1.2}). 
For simplicity, we call $p\in{\Bbb R}^2$  a {\it periodic point\/} of
$T$ of period $n$, if it is the fixed point of the iterate $T^n$. For each positive integer $n$, we put
\begin{equation}\label{eq:beta}
n\alpha=2g\pi+\beta, \quad -\pi<\beta<\pi, \ g\in{\Bbb Z}.
\end{equation}
Let $C_{n,j} $ be the circle centered at the origin with radius $r_j$ determined  by
$$
nr_j^{2s}=-\beta+2\pi j,\quad j=1,2,\ldots.
$$
It is easy to see that for $0<\beta<\pi$, the set of periodic points of $T$
with period $n$ is the disjoint union of $
C_{n,j} \ (j\geq 1).
$
For $-\pi<\beta<0$, the set of periodic points contains an
extra circle 
$$
C_{n,0}=\left\{z|n(z\overline z)^s=-\beta\right\},
$$
which is the smallest circle of periodic points of $T$ with period $n$.

A small perturbation of twist mapping may destroy the circle
$C_{n,j}$ of periodic points. However, these  circles are only deformed slightly into
Birkhoff curves, i.\ ~e.\ curves translated radially by the $n$-th iterate of the perturbation. We notice that for a perturbation of twist mapping, the Birkhoff curve can be constructed near each 
periodic circle $C_{n,j} $ which is close to the origin (see~\cite{siegelmoser}, in particular 
pages 176-177). 
However, we shall focus on the periodic points of a perturbed mapping 
near the circle $C_{n,0}$.

 Consider the complexification of the mapping (\ref{eq:1.3}), or more
generally, a holomorphic transformation of ${\Bbb C}^2$ defined by
\begin{equation}\label{eq:3.1}
\begin{split}&
\xi^\prime=e^{i\omega(\xi\eta)}\xi\bigl(1+p(\xi,\eta)\bigr),\\ &
\eta^\prime=e^{-i\omega(\xi\eta)}\eta\bigl(1+p(\xi,\eta)\bigr)^{-1},
\end{split}
\end{equation}
where $\omega(\xi\eta)=\alpha+(\xi\eta)^s$ and  $p=O(2s+1)$ is a power series   converging on 
$$\Delta_R=\bigl\{(\xi,\eta); |\xi|<R,|\eta|<R\bigr\},\quad 0<R<1.$$
Set
$$
\| p\|_R=\max_{|\xi|,|\eta|\leq R}|p(\xi,\eta)|.
$$
 For a positive number $\epsilon$, we put
$$
\Omega_{n,\epsilon}=\bigl\{(\xi,\eta)\in {\Bbb C}^2; 0<|\xi|,\ |\eta|<\epsilon_0 n^{-1/(2s)}, 1/4<|\xi/\eta|<4\bigr\}.
$$

We have the following.
\begin{thm}\label{thm:3.1}{\it
Let $\varphi$ be given by $(\ref{eq:3.1})$ with
$\| p\|_R\leq m_0.$    Then there exist positive constants $\epsilon$ and $\delta,$  depending
only on $R, m_0$ and $s,$ such that for any positive integer $n,$  the fixed
points  of $\varphi^n$ in the domain $\Omega_{n,\epsilon}$
are the union of $s$
holomorphic curves if $\beta$ in $(\ref{eq:beta})$ satisfies $-\delta<\beta<0.$
 Furthermore$,$ if $\varphi$ satisfies the reality condition
$\rho\varphi=\varphi\rho ,$ then the periodic points in the totally real
space ${\Bbb R}^2:\eta=\overline \xi$ form a closed real analytic curve$.$}
\end{thm} 
We shall see that, as a rule, the above holomorphic curves cross the level
curve $\xi\eta=c$. This is in contrast to 
the theory of area-preserving transformations,
of which 
the periodic points  are
always contained in the level curves of its integrals.

\vspace{3ex}\setcounter{section}{3}\noindent
{3.\hspace{.5em}\it Estimates for iterates}\\
\noindent
 \setcounter{thm}{0}\setcounter{equation}{0}
Let $\varphi$ be defined by (\ref{eq:3.1}). Consider the $k$-th iterate  
\begin{equation}\label{eq:3.0}\varphi^k:\begin{split}&
\xi_k=\xi e^{ik\omega(\xi\eta)}\bigl(1+p_k(\xi,\eta)\bigr)
,\\ &
\eta_k=\eta e^{-ik\omega(\xi\eta)}\bigl(1+p_k(\xi,\eta)\bigr)^{-1}
\end{split}
\end{equation}
 with $p_1(\xi,\eta)=p(\xi,\eta)$.
By setting $p_0(\xi,\eta)\equiv 0$, we have
\begin{equation}\label{eq:4.1}
p_{k+1}=p_k+p\circ\varphi^k+ p\circ\varphi^k\cdot  p_k
.
\end{equation}

We introduce some notations. For a power series $f(\xi,\eta)$ with coefficients $f_{i,j}$, let us denote
$$
\hat f(\xi,\eta)=\sum_{i,j\geq 0}|f_{i,j}|\xi^i\eta^j.
$$
We say that $g $ is a {\it majorant\/} of $f $, symbolically $f\prec g$, if 
$|f_{i,j}|\leq g_{ij}$ for all $i, j\geq 0$.
Assume that $f$ converges on $\Delta_R$. By the Cauchy inequality, one gets
$$
|f_{i,j}|\leq \frac{\| f\|_R}{R^{i+j}}.
$$
Assume further that  $f$ starts with terms of order $k$. Then one has
\begin{equation}\label{eq:4.2}
f(\xi,\eta)\prec \sum_{i+j\geq k}\frac{\| f\|_R}{R^{i+j}}\xi^i\eta^j\prec
\frac{\| f\|_R}{R^k}\cdot\frac{(\xi+\eta)^k}{1-R^{-1}\xi-R^{-1}\eta}.
\end{equation}

From (\ref{eq:4.1}), we have
$$
\hat p_{k+1}\prec \hat p_k+\hat p\circ\hat \varphi_k
+\hat p_k\cdot\hat p\circ\hat\varphi_k,
$$
where
$$\hat\varphi_k(\xi,\eta)=\frac{
 e^{k(\xi\eta)^s}}{1-\hat p_k}\left(\xi ,\eta\right).
$$
From (\ref{eq:4.2}),  we can recursively define a majorant $f_k $ for
$p_k $ by setting $f_0\equiv 0$ and 
\begin{equation}\label{eq:4.3}
f_{k+1}(\xi,\eta)= f_k(\xi,\eta)+\frac{\| p\|_R}{R^{2s+1}}\cdot
\frac{(\xi+\eta)^{2s+1}e^{k(2s+1)(\xi\eta)^s}\bigl(1-f_k(\xi,\eta)\bigr)^{-2s-2}}{1-R^{-1}(\xi+\eta)e^{k(\xi\eta)^s}\bigl(1-f_k(\xi,\eta)\bigr)^{-1}},
\end{equation}
for $k\geq 0$.

  Put
\begin{equation}\label{eq:4.4}
d_0=\min\left\{ \frac{R^{2s+1}}{2^{6s+6}m_0},\left(\frac{1}{2n}\right)^{
\frac{1}{2s}},\frac{R}{16}\right\}.
\end{equation}
Then there exists a positive constant $c_1$ depending only on $s, m_0$ and $R$
such that
\begin{equation}\label{eq:4.3a}
c_1<nd_0^{2s}<1/2.
\end{equation}
We want  to prove the following lemma.
\begin{lemma}\label{lemma:4.2}{\it
Let  
$\varphi$  be   as in {\/\it Theorem~\ref{thm:3.1}}$.$
 Then the $k$-th iterate of $\varphi$ satisfies
\begin{equation}\label{eq:4.5a}
\varphi^k:\Delta_{d_0}\to\Delta_{4d_0},\quad 1\leq k\leq n.
\end{equation}
}
\end{lemma}
\pf From (\ref{eq:4.3a}), we have
$$|e^{ik\omega(\xi\eta)}|\leq e^{nd_0^{2s}}<2,\quad (\xi,\eta)\in\Delta_{d_0}.$$
Hence
$$
|\xi_k|\leq 2|\xi|(1+|p_k(\xi,\eta)|),\quad |\eta_k|\leq
2|\eta||1+p_k(\xi,\eta)|^{-1}.
$$ 
Since $p_k\prec f_k$, it suffices to show that $ f_k$
 converges
in a neighborhood of $ \Delta_{d_0}$,  and 
\begin{equation}\label{eq:4.6}
\| f_k\|_{d_0}\leq\frac{k}{4n}.
\end{equation}

We apply the induction on $k$. For $k=0$, the assertion is trivial since $f_0=0$.
Assume that the assertion holds for some $k<n$. This implies that
$(1-f_k)^{-1}$ converges on $\Delta_{d_0}$.
From (\ref{eq:4.4}), one has
\begin{equation}\label{eq:4.7}
d_0\leq R/16.
\end{equation}
Hence, for $(\xi,\eta)\in\Delta_{d_0}$,
$$
|R^{-1}(\xi+\eta)e^{k(\xi\eta)^s}\bigl(1-f_k(\xi,\eta)\bigr)^{-1}|
<1/2.
$$
Thus the right side
of (\ref{eq:4.3}) converges on $\Delta_{d_0}$, which gives the convergence of
$f_{k+1}$ in the  same domain. Now, one can get
$$
\left\|\frac{\| p\|_R}{R^{2s+1}}\cdot
\frac{(\xi+\eta)^{2s+1}e^{k(2s+1)(\xi\eta)^s}(1-f_k)^{-2s-2}}{1-R^{-1}(\xi+\eta)e^{k(\xi\eta)^s}(1-f_k)^{-1}}\right\|_{d_0}
\leq
\frac{m_02^{6s+5}d_0^{2s+1}}{R^{2s+1}}.
$$
From (\ref{eq:4.3a}), one can 
replace 
$d_0^{2s}$ by $\frac{1}{2n}$ and rewrite the 
right side of the above inequality as
$$
\frac{m_0}{R^{2s+1}}\cdot\frac{2^{6s+4}d_0}{n},
$$
which, by (\ref{eq:4.4}), does not exceed $\frac{1}{4n}$. This proves (\ref{eq:4.6}). 

Under assumptions as in Lemma~\ref{lemma:4.2}, we now give
some estimates for $p_k $  and $f_k $.
Notice that $p_k $ and $f_k $ start with terms of order
$2s+1$. By (\ref{eq:4.6}), the Schwarz lemma gives
\begin{equation}\label{eq:4.7c}
\mid p_k(\xi,\eta)\mid_r,\quad\mid f_k(\xi,\eta)\mid_r\leq \frac{1}{4}\bigl(\frac{r}{d_0}\bigr)^{2s+1} ,\quad  r\leq d_0.
\end{equation}

We are ready to prove Theorem~\ref{thm:3.1}.

{\bf Proof of Theorem~\ref{thm:3.1}.} Let $\varphi, m_0, d_k$ and $n$ be as  in Lemma~\ref{lemma:4.2}.
Set
$$
\xi=\zeta w,\quad \eta=\zeta w^{-1}.
$$
We require  that
\begin{equation}\label{eq:4.7a}
0<|\zeta|<d_0/2, \quad 1/2<|w|<2.
\end{equation}

Consider the equation $\xi_n=\xi$, or equivalently,
\begin{equation}\label{eq:4.8}
e^{i(\beta+n\zeta^{2s})}\{1+p_n(\zeta
w,\zeta
w^{-1} )\}=1,
\end{equation}
in which $\beta$ is  determined by (\ref{eq:beta}) for a given $n$.

From (\ref{eq:4.7c}) and (\ref{eq:4.7a}), it follows that
\begin{equation}\label{eq:4.11}
|p_n(\zeta w,\zeta w^{-1} )|\leq
1/4.
\end{equation}
We also have $n|\zeta|^{2s}<1/2$ and $|\beta|<\pi$. Thus (\ref{eq:4.8}) is
reduced to
$$
\beta+n\zeta^{2s}-i\log\bigl(1+p_n(\zeta w,\zeta w^{-1} )\bigr)=0,
$$
in which the logarithm assumes principal values. We further rewrite the above equation as
\begin{gather}\label{eq:4.12}
\zeta^{2s}\bigl(1+h(\zeta,w)\bigr)=-\frac{\beta}{n}\\
\intertext{with}
h(\zeta,w)=\frac{1}{i\,n\zeta^{2s}}
\log\bigl(1+p_n(\zeta w,\zeta w^{-1})\bigr).
\label{eq:4.13}
\end{gather}
Using (\ref{eq:4.7c}), we have
$$
|h(\zeta,w)|\leq
\frac{1}{n|\zeta|^{2s}}\left(\frac{2|\zeta|}{d_0}\right)^{2s+1}.
$$
From (\ref{eq:4.3a}), it follows that
\begin{equation}\label{eq:4.14}
|h(\zeta,w)|
\leq \frac{n^{\frac{1}{2s}}}{c_2}|\zeta|,
\end{equation}
in which $c_2<c_1$ is a  constant  depending only on $m_0,R$ and $s$.

We now  take
\begin{equation}\label{eq:4.15}
\epsilon_0=c_2
,\quad \delta=\left(\frac{c_2}{4}\right)^{2s}.
\end{equation}
Put
$$
r_0=\frac{1}{2}\epsilon_0 n^{-\frac{1}{2s}}<\frac{d_0}{2}.
$$
It follows from (\ref{eq:4.14}) and (\ref{eq:4.15}) that
\begin{equation}\label{ade-3}
|h(\zeta,w)|\leq 1/4,\quad |\zeta|<r_0.
\end{equation}
Now (\ref{eq:4.12}) is reduced to
\begin{equation}\label{eq:4.16}
\zeta\bigl(1+h(\zeta,w)\bigr)^{\frac{1}{2s}}=e^{i\frac{j\pi}{s}}
\left(-\frac{\beta}{n}\right)^{\frac{1}{2s}},\quad j=1,2,\ldots,2s.
\end{equation}
Furthermore, for $|\zeta|=r_0$, we have
\begin{equation}\label{ade-2}
|(1+h)^{\frac{1}{2s}}-1|=|\frac{h}{2s}\int_0^1(1+th)^{\frac{1}{2s}-1}\, dt|<1/6.
\end{equation}
We now assume that $|\beta|<\delta$. Then (\ref{eq:4.15}) gives
\begin{equation}\label{ade-1}
\bigl |\bigl(\frac{-\beta}{n}\bigr)^{\frac{1}{2s}}\bigr |<\frac{c_2}{4}n^{-\frac{1}{2s}}=r_0/2.
\end{equation}
 Thus, the Rouch\'e theorem implies that for each $j$,
(\ref{eq:4.16}) has a
unique solution $\zeta=\zeta_j(w)$ in the  disk $|\zeta|< r_0.$
Clearly, the solution is  holomorphic for $1/2<|w|<2$. 

Notice that the transformation $(\zeta,w)\rightarrow(\zeta w,\zeta w^{-1})$
is two-to-one. Also, the equation (\ref{eq:4.8}) is invariant for the
transformation $(\zeta, w)\rightarrow(-\zeta,-w)$. Hence, the $2s$
holomorphic curves $\zeta_j(w)\ (1\leq j\leq 2s)$ in $(\zeta,w)$--coordinates
give us $s$
holomorphic curves  $(\xi,\eta)=(\zeta_j w,\zeta_j w^{-1})\ (1\leq j\leq
s)$.
Since 
$\kappa(\xi,\eta)=\xi\eta$ is 
invariant under $\varphi$, the solutions $(\xi,\eta)$ to
$\xi_n=\xi$ give the fixed points of $\varphi^n$.  

We now assume that $\varphi$ satisfies the reality condition, and want to show that $\zeta(w)\equiv\zeta_{2s}(w)$ is a real
valued function
for $|w|=1$.
Since $\rho\varphi=\varphi\rho$, we have $\varphi^n\rho=\rho\varphi^n$.
Hence
$$1+\overline p_n(\eta,\xi)=\bigl(1+p_n(\xi,\eta)\bigr)^{-1}.
$$
It is easy to see that for  $|w|=1$,
$$
\overline h(\overline \zeta, \overline w)=h(\overline
\zeta,w ).
$$
Conjugating (\ref{eq:4.16}), we get
$$
\overline \zeta\bigl(1+h(\overline
\zeta,w )\bigr)^{\frac{1}{2s}}=e^{-i\frac{j\pi}{s}}\left(-\frac{\beta}{n}\right)^{ \frac{1}{2s}},$$
which is precisely the equation (\ref{eq:4.16}) for $j=s$. From the uniqueness of the solution, we obtain that
$\zeta(w)=\overline{\zeta(w )}$ for $|w|=1$.
The proof of
Theorem~\ref{thm:3.1} is complete.

\vspace{3ex}\setcounter{section}{4}\noindent
{4.\hspace{.5em}\it Analytic dependence on coefficients}\\
\noindent
 \setcounter{thm}{0}\setcounter{equation}{0}
We shall first discuss the analytic dependence of holomorphic curves of
periodic points of $\varphi$ defined
by (\ref{eq:3.1}) on the coefficients of 
$p$. We need  power series in infinitely many variables. The
convergence of such a power series will always mean that it converges absolutely.

For each positive $m_0$ and $R$, we define 
$$
D^\infty(m_0,R)=\bigl\{p=(p_{i,j})  ;|p_{i,j}|\leq\frac{m_0}{4(2R)^{i+j}},i+j>2s\bigr\} .
$$
For each $p\in D^\infty(m_0,R)$, we put
$$
\widetilde p(\xi,\eta)=\sum_{i+j>2s}p_{i,j}\xi^i\eta^j.
$$
Then  $\widetilde p(\xi,\eta)$   converges on
$\Delta_R$ and 
$$
\|\tilde p\|_R\leq \sum_{i+j>2s}|p_{i,j}|R^{i+j}
\leq\frac{m_0}{4}\sum_{k>2s}(k+1)\frac{1}{2^k}<m_0.
$$

Consider the transformation $\varphi$ defined by (\ref{eq:3.1}) with $p(\xi,\eta)=\tilde p(\xi,\eta)$. Put the iterate $\varphi^k$ in the form (\ref{eq:3.0}) with $p_k(\xi,\eta)=
p_k(\xi,\eta,p)$.
 We also denote the solution $\zeta$ to (\ref{eq:4.16})
by $\zeta_j(w,p )$.
We have seen that for a fixed $p\in D^\infty(m_0,R)$,  $p_k(\xi,\eta,p)$ 
 converges in $\xi$ and $\eta$, and $\zeta(w ,p)$ converges as a
Laurent series in $w$.
 Next, we want to show that they actually
converge  as series expansions with respect to all variables. For $p\in D^\infty(m_0,R)$, we put
$p_0^*(\xi,\eta,p)\equiv 0$ and define recursively
$$
p_{k+1}^*=p_k^*+\widetilde p\circ\varphi_k^*(1+p_k^*),\quad k=0,1,\ldots, n-1,
$$
in which
$$
\varphi_k^*(\xi,\eta)=\bigl(\xi e^{k(\xi\eta)^s}(1-p_k^*)^{-1},\eta
e^{k(\xi\eta)^s}(1-p_k^*)^{-1}\bigr).
$$
From (\ref{eq:4.1}), it follows that
$$
p_k(\xi,\eta,p)\prec p_k^*(\xi,\eta,p).
$$
Put 
$$p^0(\xi,\eta)=\sum_{i+j>2s}
\frac{1}{4}m_0(2R)^{-i-j}\xi^i\eta^j.$$
Clearly, for a fixed $p\in D^\infty(m_0,R)$, one has
$$
p_k^*(\xi,\eta,p)\prec p_k^*(\xi,\eta,p^0)\prec f_k^*(\xi,\eta).
$$
in which $f_k^*(\xi,\eta)$ are defined through the recursive formulas (\ref{eq:4.3}) by setting 
$p=p^0$ and $f_0^*\equiv 0$. From (\ref{eq:4.6}), it follows that
  $p_k^*(\xi,\eta,p^0)$
converges on $\Delta_{d_0}$ with $\|p_k^*(\cdot,\cdot,p^0)\|_{d_0}\leq 1/4$. Hence, $p_k^*(\xi,\eta,p)$ 
  converges on $\Delta_{d_0}\times D^\infty(m_0,R)$ with
\begin{equation}\label{ade1}
\|p_k^*(\cdot,\cdot,p)\|_{d_0}< 1/4,\qquad \mbox{for}\ p\in D^\infty(m_0,R).
\end{equation}

To show the convergence of $\zeta_j(w, p)$, we use the Residue formula
$$
\zeta_j(w,p)=\frac{1}{2\pi i}\int_{|\zeta|=r_0}\frac{\partial_\zeta F_j(\zeta,w,p)}{F_j(\zeta,w,p)}d\zeta,
$$
where
$$
F_j(\zeta,w,p)=\zeta\left(1+h(\zeta,w,p)\right)^{\frac{1}{2s}}-e^{i\frac{j\pi}{s}}\left(-\frac{\beta}{n}\right)^{\frac{1}{2s}},
$$
and $h(\zeta,w,p)$ is given by (\ref{eq:4.13}). As series expansion in $\zeta, w$ and $p$, we have
$$
h\prec\frac{p_n^*}{n\zeta^{2s}(1-p_n^*)}\equiv \widetilde h.
$$
As (\ref{ade-3}), one can use (\ref{ade1}) to get
\begin{equation}\label{adeh}
|p_n^*(r_0,w,p)|<1/4,\quad \mbox{for}\ 1/2<|w|<2,\ |\zeta|\leq r_0,\ p\in D^\infty(m_0,R).
\end{equation}
From (\ref{ade-2}) and (\ref{ade-1}),  it follows that  for each $p\in D^\infty(m_0,R)$, the function $1/F_j$ has the expansion
$$
\frac{1}{\zeta}\sum_{k=0}^\infty\bigl((1+h(\zeta,w,p))^{\frac{1}{2s}}-1+
\frac{e^{i\frac{j\pi}{s}}}{\zeta}\bigl(-\frac{\beta}{n}\bigr)^{\frac{1}{2s}}\bigr)^k,\quad |\zeta|=r_0,\ 1/2<|w|<2.
$$
Since $(1+x)^{\frac{1}{2s}}\prec\frac{1}{1-x}$, then 
$$
(1+h)^{\frac{1}{2s}}\prec \frac{1}{1-\widetilde h}.
$$
Hence, we see that as series expansion in $\zeta, w, p$, $1/F_j$ is majorized by
\begin{equation}
\label{ade2}
\frac{1}{\zeta}\sum_{k=0}^\infty\left(\frac{\widetilde h(\zeta,w,p)}{1-\widetilde h(\zeta,w,p)} +
\frac{ 1}{ \zeta}\left(\frac{|\beta|}{n}\right)^{\frac{1}{2s}}\right)^k.
\end{equation}
From (\ref{ade-1}) and (\ref{adeh}), it follows that $1/F_j$ converges 
  for $1/2<|w|<2, p\in D^\infty(m_0,R)$ and $|\zeta|=r_0$. By the Residue formula, we 
see that $\zeta_j(w,p)$ converges 
 for $1/2<|w|<2$ and $ p\in D^\infty(m_0,R)$.
  We have proved the following. 
\begin{lemma}\label{lemma:5.1}{\it
The solution $\zeta_j(w, p)$ to $(\ref{eq:4.16})$ converges for
$1/2<|w|<2$ and $p\in D^\infty(m_0,R).$
}
\end{lemma}

We now want to apply Lemma~\ref{lemma:5.1} to the family of reversible
transformations (\ref{eq:1.3}).
Let us denote $$D^\infty=\{a=(a_{i,j})\mid i+j>2s ,|a_{i,j}|\leq 1\}.$$
Fixing $a\in D^\infty$, we consider the holomorphic mapping 
\begin{equation}\label{eq:5.1}
\phi_a\colon\begin{split}&
\xi^\prime=e^{i\tilde a(\xi,\eta)}\xi,\\ &
\eta^\prime=e^{-i\tilde a(\xi,\eta)}\eta,
\end{split}
\end{equation}
with $\tilde a(\xi,\eta)=\sum_{i+j>2s}a_{i,j}\xi^i\eta^j.$
Let us denote
\begin{equation*}
\phi_a^*\colon\begin{split}&
\xi^\prime=e^{\tilde a(\xi,\eta)}\xi,\\ &
\eta^\prime=e^{\tilde a(\xi,\eta)}\eta.
\end{split}
\end{equation*}
We can put 
\begin{equation}\label{eq:5.2}
\phi_a^{-1}\colon\begin{split}&
\xi^\prime=e^{-ib(\xi,\eta)}\xi,\\ &
\eta^\prime=e^{ib(\xi,\eta)}\eta,
\end{split}
\end{equation}
for some formal power series $b(\xi,\eta)$ without the constant 
term. This leads to the following identity
\begin{equation}\label{eq:5.3}
b(\xi,\eta)=\tilde a\bigl(e^{-ib(\xi,\eta)}\xi,e^{ib(\xi,\eta)}\eta\bigr).
\end{equation}
We consider $b(\xi,\eta)$ as a power series in $\xi,\eta$ and 
$a$. Then  $b(\xi,\eta)$ is majorized by
\begin{equation}\label{eq:eq:5.4}
 b^*(\xi,\eta,a)= 
\tilde a\bigl(e^{ b^*(\xi,\eta,a)}\xi,e^{ b^*(\xi,\eta,a)}\eta\bigr).
\end{equation}
Put
\begin{equation*}
\psi_a\colon\begin{split}&
\xi^\prime=e^{b^*(\xi,\eta,a)}\xi,\\ &
\eta^\prime=e^{b^*(\xi,\eta,a)}\eta.
\end{split}
\end{equation*}

We have the following.
\begin{lemma}\label{lemma:5.2}{\it
Let $\psi_a(\xi,\eta)$ be given as   above$.$ Then there exists a constant $R_1,0<R_1<1,$ 
independent of $a\in D^\infty,$ such that  $\psi_a $  is given by power series which are convergent
 for $(\xi,\eta,a)\in\Delta_{R_1}\times D^\infty.$
 Moreover$,$ for each $a\in D^\infty,$  
\begin{equation}\label{eq:5.5}
\psi_a\colon\Delta_{R_1}\to\Delta_{1/3}.
\end{equation}
}
\end{lemma}
\pf
Put 
$$
e(\xi,\eta)=\sum_{i+j\geq 1}\xi^i\eta^j.
$$
By the implicit function theorem, there exists a positive 
number $R_1<1/6$ such that
$b^*(\xi,\eta,e)$ converges on $\Delta_{R_1}$, and 
$|b^*(\xi,\eta,e)|<1/2$. Clearly, we have (\ref{eq:5.5}).
 Now for each $a\in D^\infty$, we have $a(\xi,\eta)\prec e(\xi,\eta)$.
Since $b^*(\xi,\eta,a)$ has positive coefficients, 
this shows that $b^*(\xi,\eta,a)$ converges on 
$\Delta_{R_1}\times D^\infty$. The proof of
the lemma is complete.

We are ready to prove the following.
\begin{prop}\label{prop:5.3}{\it
Let $\varphi_a$ be  given by $(\ref{eq:1.3}).$ Then there
exist positive constants $\epsilon_0$ and $ \delta$ which are independent of $n$ and $a $ such that in the
polar coordinates $(r,\theta),$  the fixed points of  $\varphi_a^n$ in $|z|<\epsilon_0 n^{-1/(2s)}$ form a closed analytic 
curve $r=\zeta(w, a),|w|=1,$ provided 
$\beta$  in $(\ref{eq:beta})$ satisfies $-\delta<\beta<0.$ Furthermore$,$ $\zeta(w ,a)$ converges for  $1/2<|w|<2$ and $a\in\sum.$
}
\end{prop}
\pf For each $a\in\sum$, let $b(\xi,\eta)$ be given by (\ref{eq:5.3}). We have  $b(z,\overline z)=\overline b(\overline z, z)$. Hence 
$$
\varphi_a=\phi_a\circ T\circ\overline\phi_a^{-1}=\phi_a\circ T\circ\phi_b.$$
Put
$$
\phi_a^*(\xi,\eta)=\bigl(\xi e^{\tilde a(\xi,\eta)},\eta e^{\tilde a(\xi,\eta)}\bigr),
\quad  T^*(\xi,\eta)=\bigl(\xi e^{(\xi\eta)^s},\eta e^{(\xi\eta)^s }\bigr).
$$
It is easy to see that as power series in $\xi=z,\eta=\overline z$ and $a$, one has
$$
\varphi_a(\xi,\eta)\prec\phi_a^*\circ T^*\circ\psi_a(\xi,\eta)=\varphi_a^*(\xi,\eta).
$$
For a fixed $a\in\sum$,  $\varphi_a^*$ is majorized by
$\varphi_e^* $. Since $\varphi_e^*$ is
convergent in the disk $\Delta_{R_1}$, then  $\varphi_a(z,\overline z)$
converges on $\Delta_{R_1}\times D^\infty$.
  Put $\varphi_a$ in the form (\ref{eq:3.1}) with   $p(\xi,\eta,a)$ in place of $p(\xi,\eta)$.   Then $p(\xi,\eta,a)$  converges in
$\Delta_{R_1}\times D^\infty$. Therefore, there exists a constant  $m_0$, say $m_0=4\hat p(R_1,R_1,e)$, such that for each $a\in\sum$
$$
p(\xi,\eta,a)\in D^\infty(R_1/2,m_0) .
$$
Now Theorem~\ref{thm:3.1} and Lemma~\ref{lemma:5.1} give us the required convergence of $\zeta_j(w,a)$. The proof of  Proposition~\ref{prop:5.3} is complete.

\vspace{3ex}\setcounter{section}{5}\noindent
{5.\hspace{.5em}\it Proof of  Theorem\/~$\ref{thm:1.1}$}\\
\noindent
\setcounter{thm}{0}\setcounter{equation}{0}
  From (\ref{eq:4.1}),  one
has  
\begin{equation}\label{eq:6.1}
p_n(\xi,\eta,p)=\sum_{k=0}^{n-1}p\bigl(\xi e^{i\,k\omega(\xi\eta)},\eta e^{-i\,
k\omega(\xi\eta)}\bigr)+O(p^2),
\end{equation}
in which $ O(p^2)$ stands for terms  of order at least $2$ in variables
$p_{i,j}$.  Now (\ref{eq:4.13}) gives 
\begin{equation}\label{eq:6.2}
h(\zeta,w,p)=\frac{1}{i\, n\zeta^{2s}}\sum_{k=0}^{n-1}p\bigl(\zeta w
e^{i\,k\omega(\zeta^2)},\zeta w^{-1}e^{-i\,
k\omega(\zeta^2)}\bigr)+O(p^2).
\end{equation}

Let $\zeta_n(w, p)$ be the solution to the equation
(\ref{eq:4.16}) with $j=2s$.
Denote by $K\zeta(w, p)$ the  constant term of $\zeta(w,p)$  with respect to
variables
$p_{i,j}$, and by $L\zeta(w, p)$   the   linear part of $\zeta(w,p)$.  Then from (\ref{eq:4.16}), we get
\begin{equation}\label{eq:6.3}
K\zeta(w)=\zeta_0,\quad \zeta_0=\left(-\frac{\beta}{n}\right)^{\frac{1}{2s}}.
\end{equation}
Now (\ref{eq:6.2})  gives
\begin{gather}\label{eq:6.4}
L\zeta(w, p)=\frac{i}{ 2sn\zeta_0^{2s}}\sum_{k=0}^{n-1}p(\zeta_0 wu^k,\zeta_0 w^{-1}u^k)\\
\intertext{with}
u=e^{i\omega(\zeta_0^2)}.\nonumber
\end{gather}

We now can complete the proof of Theorem~\ref{thm:1.1}. 
 For each $a\in\sum$, write
\begin{equation*}
\varphi_a=\phi_a\circ T\circ\overline\phi_a^{-1}\colon
\begin{split}&
\xi^\prime=\xi e^{i\omega(\xi\eta)}\bigl(1+p(\xi,\eta,a)\bigr),\\ &
\eta^\prime=\eta e^{-i\omega(\xi\eta)}\bigl(1+p(\xi,\eta,a)\bigr)^{-1}.
\end{split}
\end{equation*}
From (\ref{eq:5.2}) and (\ref{eq:5.3}), we get
$$
\phi_a^{-1}(\xi,\eta)=\bigl(\xi-i\xi \tilde a(\xi,\eta),\eta+i\eta \tilde a(\xi,\eta)\bigr)+O(a^2),
$$
in which  $O(a^2)$ terms of order at least $2$ in variables
 $a_{i,j}$.
From $\tilde a(\xi,\eta)=\overline{\tilde a}(\eta,\xi)$,
 We have
$$
\overline \phi_a^{-1}(\xi,\eta)=\bigl(\xi+i\xi\tilde a(\eta,\xi),\eta-i\eta \tilde a(\eta,\xi)\bigr)+O(a^2).
$$
Now one   gets
\begin{equation}\label{eq:6.5}
p(\eta,\xi,a)=i\tilde a\bigl(\xi e^{i\omega(\xi\eta)},\eta
e^{-i\omega(\xi\eta)}\bigr)+i\tilde a(\eta,\xi)+O(a^2).
\end{equation}

Let $\epsilon_0,\delta$ be as in Proposition~\ref{prop:5.3}. Choose
a sequenece of positive integers $n_k\to\infty$ such that for some
$g_k\in{\Bbb Z}$,
$$
n_k\alpha=2g_k\pi+\beta_k,\quad -\delta<\beta_k<0.
$$
Then for each   $a\in\sum$, the periodic points of
$\varphi_a$ in $|z|<\epsilon_0 n_k^{-1/(2s)}$ 
form an analytic
curve $r=\zeta_{n_k}(w, a)$, $|w|= 1$,    given by
$$
\zeta_{n_k}\left(1+h(\zeta_{n_k},w, a)\right)^{\frac{1}{2s}}=\left(-\frac{\beta_k}{n_k}
\right)^{\frac{1}{2s}}\equiv\zeta_0.
$$
By Proposition~\ref{prop:5.3}, $\zeta_{n_k}(w, a)$ converges for
$1/2<|w|<2$ and $a\in\sum$. From
 (\ref{eq:6.4}) and (\ref{eq:6.5}), one sees that with respect to the variables $a_{i,j}$, the linear part of $\zeta_{n_k}$ is
$$
L\zeta_{n_k}(w, a)=-\frac{1}{ 2sn_k\zeta_0^{2s}}\sum_{j=0}^{n_k-1}\tilde a(\zeta_0 w u^{j+1},\zeta_0 w^{-1}\overline  u^{j+1})
+\tilde a(\zeta_0 w^{-1}
\overline u^j,\zeta_0 w u^j).
$$
Let $H_k(a)$ be the coefficients of $w^{n_k}$ in  the Laurent expansion
 of
$\zeta_{n_k}(w, a)$ with respect to 
$w$. Notice that $u^{n_k}=1$. Then we have
$$
H_k(a)=-\zeta_0^{n_k-2s}
(a_{n_k,0}+a_{0,n_k})/(2s)+h_k(a),
$$
in which $h_k(a)$ contains no linear terms in $a_{j,0}$ and
$a_{0,j}$ for all $j$.
It is clear that $H_k(a)\ (k=1,2,\ldots)$ are functionally independent.

 Fix $a\in\sum$. Assume that  $H_k(a)\neq 0$. 
Since $H_k(a)$ is the 
coefficient of $w^{n_k}$ in the expansion of  $\zeta_{n_k}(w, a)$,
  then  $\zeta_{n_k}(w, a)$ is not constant with respect to the variable $w$.  Let $I_{k}=\{\zeta_{n_k}(w, a);|w|=1\}$.  Then $I_{k}$ is a closed interval of positive length. For each $r\in I_{k}$, $\varphi_a ^{n_k}$ has a  fixed point on $|z|=r$. Obviously,  $ I_{k}$  is contained in any given interval $(0,\epsilon)$ for large $k$. The proof of Theorem~\ref{thm:1.1} is complete.

\vspace{3ex}\setcounter{section}{6}\noindent
{6.\hspace{.5em}\it Complex tangents of real surfaces}\setcounter{thm}{0}\setcounter{equation}{0}\\
\noindent
It is natural to replace the reality condition of a reversible
transformation of ${\Bbb C}^2$ by the reversibility of the transformation
with respect to an anti-holomorphic involution. This means that one may 
ask for the classification of biholomorphic mappings which are reversible
through  a holomorphic involution  as well as an anti-holomorphic involution.
Such reversible holomorphic mappings arise from the real analytic surfaces
with non-degenerate complex tangents, which is the main subject to be
discussed in this section.

A real analytic surface $M$ has a {\it complex tangent\/} at $p$ if $T_pM$ is a complex line
in ${\Bbb C}^2$. Introduce local coordinates such that $p=0$,  and 
 $M:z_2=az_1^2+bz_1\overline z_1+c\overline z_1^2+O(|z|^3)$. The complex tangent is said to be  
  {\/ \it non-degenerate\/}  if  $b\neq 0$. Then
for a further change of local holomorphic coordinates near $ 0$, one may
assume that M is given by (\ref{eq:1.4}). The $\gamma$ in (\ref{eq:1.4}) is the Bishop invariant~\cite{bishop}. The
complex tangent is said to be {\it elliptic$,$ parabolic\/} or {\/ it
hyperbolic\/} if $0\leq \gamma<1/2, \gamma=1/2$ or $ \gamma>1/2$,
respectively.

To describe a pair of intrinsic involutions introduced in~\cite{moserwebster}, we consider the
complexified surface in ${\Bbb C}^2\times{\Bbb C}^2$
\begin{equation*}
M^c:\begin{split}&
z_2=z_1w_1+\gamma z_1^2+\gamma w_1^2+q(z_1,w_1),\\ &
w_2=z_1w_1+\gamma z_1^2+\gamma w_1^2+\overline q(w_1,z_1).
\end{split}
\end{equation*}
Let $\pi_1(z,w)=w$ and $\pi_2(z,w)=z$ with $z=(z_1,z_2), w=(w_1,w_2)$.
 Then for
$\gamma\neq 0$, $\pi_j|_{M^c}$ is a branched
double-covering. The covering transformation $\pi_j|_{M^c}$ gives an
involution $\tau_j$. Two involutions $\tau_1$ and $\tau_2$ are conjugate to
each other by the anti-holomorphic
involution $(z,w)\mapsto(\overline w,\overline z)$ restricted to $M^c$. After a
change  of local holomorphic coordinates of $M^c\equiv {\Bbb C}^2$, we may
assume that $\tau_j\ (j=1,2)$ is  given by (\ref{eq:2.15}). 
   
Let us restrict  
ourselves to the hyperbolic case. Then we have
$\lambda_1=\lambda=\lambda_2^{-1}$. The number $\lambda$ and  the
Bishop invariant $\gamma$ are related by
$$
\gamma\lambda^2-\lambda+\gamma=0.
$$
From $\gamma>1/2$, it follows that $|\lambda|=1$. We say that $\gamma$ is {\it exceptional\/} if
$\lambda $ is a root of unity. Furthermore, the anti-holomorphic involution
$(z,w)\mapsto  (\overline w,\overline z)$ restricted to $M^c$
is given by
\begin{equation}\label{eq:7.1}
\rho(\xi,\eta)=\bigl(\overline \xi,\overline \eta\bigr).
\end{equation}
Thus, the reality condition  on $\{\tau_1,\tau_2\}$ is given by
\begin{equation}\label{eq:7.2}
\rho\tau_1=\tau_2\rho.
\end{equation}

The importance of the triple $\{\tau_1,\tau_2,\rho\}$, described above, is that it
characterizes a real analytic surface as follows. A triple
$\{\tau_1,\tau_2,\rho\}$,  given by (\ref{eq:2.15}), (\ref{eq:7.1}) and (\ref{eq:7.2}), always generates a real analytic
surface (\ref{eq:1.4}).
Moreover, two real analytic surfaces are biholomorphically equivalent, if
and only if their corresponding pairs of involutions are equivalent
through a biholomorphic transformation $\Phi$ satisfying the reality
condition 
$$
\Phi\rho=\rho\Phi.
$$
We refer to~\cite{moserwebster} for details.

For $a\in D^\infty$, we define an involution
\begin{equation}\label{eq:7.3}
\tau_1=\phi_a\circ T_1\circ\phi_a^{-1},
\end{equation}
in which $\phi_a$ is given by (\ref{eq:5.1}), and
$$
T_1(\xi,\eta)=\bigl(e^{\frac{1}{2}i\omega(\xi\eta)}\eta,
e^{-\frac{1}{2}i\omega(\xi\eta)}\xi\bigr).
$$
Consider
\begin{equation}\label{eq:7.4}
\varphi\equiv\tau_1\tau_2,\quad \tau_2=\rho\tau_1\rho,\quad T_2=\rho T_1\rho.
\end{equation}
Then as power series in $\xi,\eta, a$ and $\overline a$, one has
$$
\varphi(\xi,\eta)\prec \phi_a^*T_1^*\psi_a^*\overline \phi_a^*T_2^*\overline \psi_a^*(\xi,\eta)\equiv \psi(\xi,\eta).
$$
By Lemma~\ref{lemma:5.2},  we know that there is $0<R_2<1$ such that
$\varphi(\xi,\eta)$ converges for $(\xi,\eta,a ,\overline a)\in\Delta_{R_2}\times D^\infty\times\overline D^\infty $. 
Put $\varphi$  in the form (\ref{eq:3.1}).
Thus $p(\xi,\eta)=p(\xi,\eta,a,\overline a)$ converges on
$\Delta_{R_1}\times D^\infty\times \overline D^\infty $. By
Lemma~\ref{lemma:5.1}, there are positive constants $\epsilon$ and
$\delta_0$ such that if $\beta$ in (\ref{eq:beta}) satisfies
$-\delta_0<\beta<0$, then the fixed points of $\varphi^n$ in the domain 
$
\Omega_{n,\epsilon} $
 are the union of $s$ holomorphic curves
 \begin{equation}\label{eq:7.5}
\gamma_{n,j}\colon \xi=\zeta_jw,\quad \eta=\zeta_jw^{-1},\qquad \zeta_j=\zeta_j(w, a,\overline a),\quad 1\leq j\leq s.
\end{equation}
Moreover, $\zeta=\zeta_j(w, a,\overline a)$ converges for $1/2<|w|<2$
and $(a,\overline a)\in D^\infty\times \overline D^\infty $. 

The following result implies that in general the periodic points 
 are not contained in the totally real space
$${\Bbb R}^2\colon\overline \xi=\xi,\qquad \overline \eta=\eta . 
$$
 We have
\begin{thm}\label{thm:7.1}{\it
Let  $\epsilon,\delta_0,n$ and $\Omega_{n,\epsilon}$ be as above$.$ Then there is a 
power series $H(a,\overline a)$ converging on $D^\infty\times
\overline D^\infty$ such that  $\gamma_{n,j}$ intersects the totally real space ${\Bbb R}^2$
 at isolated points$,$
provided   $H(a,\overline a)\neq 0.$
}
\end{thm}

For the proof of Theorem~\ref{thm:7.1},  we need the second order power expansion of
$\zeta_j(w, a,\overline a)$ with respect to variables $a$ and
$\overline a$. In the following discussion, we shall denote  $Lp$ the linear
terms
of  a power series $p$ in $a$ and
$\overline a$,
 and $Qp$ the second order terms of  $p$.

 To simplify the computation, we first assume that 
$$
\widetilde a(\xi,\eta)=a_{n,0}\xi^n. 
$$
Let $\phi_a$ be defined by (\ref{eq:5.1}). Then one has
\begin{equation}\label{eq:7.6}
\phi_a:
\begin{split}&
\xi^\prime=\xi\bigl(1+ia_{n,0}\xi^n-\frac{1}{2}a_{n,0}^2\xi^{2n}+
O(a_{n,0}^3)\bigr),\\ &
\eta^\prime=\eta\bigl(1-ia_{n,0}\xi^n-\frac{1}{2}a_{n,0}^2\xi^{2n}+
O(a_{n,0}^3)\bigr).
\end{split}
\end{equation}
Denote the inverse of $\phi_a(\xi,\eta)$ by $\phi_{-b}(\xi,\eta)$. Clearly,
we have $b(\xi,\eta)=a_{n,0}\xi^n+O(a_{n,0}^2)$. From (\ref{eq:5.3}), it follows that
$$
b(\xi,\eta)=a_{n,0}\xi^n-ina_{n,0}^2\xi^{2n}+O(a_{n,0}^3).
$$
Thus, one has
\begin{equation}\label{eq:7.7}
\phi_a^{-1}:
\begin{split}&
\xi^\prime=\xi\bigl(1-ia_{n,0}\xi^n-\frac{2n+1}{2}a_{n,0}^2\xi^{2n}+
O(a_{n,0}^3)\bigr),\\ &
\eta^\prime=\eta\bigl(1+ia_{n,0}\xi^n+\frac{2n-1}{2}a_{n,0}^2\xi^{2n}+
O(a_{n,0}^3)\bigr).
\end{split}
\end{equation}

Noticing that $\kappa(\xi,\eta)=\xi\eta$ is invariant under $\tau_1$, we can put
\begin{equation*}
\tau_1:
\begin{split}&
\xi^\prime=e^{\frac{1}{2}i\omega}\eta\bigl(1+q(\xi,\eta)\bigr),\\ &
\eta^\prime=e^{-\frac{1}{2}i\omega}\xi\bigl(1+q(\xi,\eta)\bigr)^{-1}.
\end{split}
\end{equation*}
From (\ref{eq:7.3}), (\ref{eq:7.6}) and (\ref{eq:7.7}), it follows that 
\begin{equation*}
\begin{split}
q(\xi,\eta)=
&ia_{n,0}(\xi^n+e^{\frac{1}{2}in\omega}\eta^n)
-\frac{1}{2}a_{n,0}^2(\xi^n+e^{\frac{1}{2}in\omega}\eta^n)^2
\\
& +na_{n,0}^2\xi^n(\xi^n-e^{\frac{1}{2}in\omega}\eta^n)
+O(a_{n,0}^3).\end{split}
\end{equation*}
From
(\ref{eq:7.4}), we get
\begin{equation*}
\tau_2:
\begin{split}&
\xi^\prime=e^{-\frac{1}{2}i\omega}\eta\bigl(1+\overline q(\xi,\eta)\bigr),\\ &
\eta^\prime=e^{\frac{1}{2}i\omega}\xi\bigl(1+\overline q(\xi,\eta)\bigr)^{-1}.
\end{split}
\end{equation*}
Put $\varphi=\tau_1\tau_2$ in the form (\ref{eq:3.1}) with 
$$
p(\xi,\eta)=\bigl(1+\overline
q(\xi,\eta)\bigr)^{-1}\bigl(1+q\circ\tau_2(\xi,\eta)\bigr)-1.
$$
Clearly, one has
\begin{equation}\label{eq:7.8}
Lp(\xi,\eta)
=i(e^{in\omega}a_{n,0}+\overline a_{n,0})\xi^n+ie^{-\frac{n}{2}i\omega}(a_{n,0}+\overline a_{n,0})\eta^n.
\end{equation}
 Notice that
$$
Q(q\circ\tau_2)=Qq\circ T_2+L\overline q\left(\eta e^{-\frac{1}{2}i\omega}
\partial_\xi Lq\circ T_2-\xi e^{\frac{1}{2}i\omega}\partial_\eta Lq\circ T_2
\right).
$$
 Then we get
\begin{equation}\label{eq:7.9}
\begin{split}Qp=
 & Qq\circ T_2
-Lq\circ T_2\cdot L\overline q-Q\overline q+(L\overline q)^2\\ 
& +L\overline q\left(\eta e^{-\frac{1}{2}i\omega}
\partial_\xi Lq\circ T_2-\xi e^{\frac{1}{2}i\omega}\partial_\eta Lq\circ T_2
\right).
\end{split}
\end{equation}

 We now assume that
$
4s|n.
$
Let $\zeta_j(0)$ be the right side of (\ref{eq:4.16}). Then from (\ref{eq:beta}), it is easy to see that
$$
e^{\frac{i}{2}n\omega(\zeta_j(0)^2)}=1.
$$
 Thus (\ref{eq:7.9}) gives
\begin{equation}\label{eq:7.10}
Qp(\zeta_j(0) w,\zeta_j(0)w^{-1})=
n a_{n,0}^2\bigl(\zeta_j(0)w\bigr)^{2n}+q_1(a_{n,0},\overline a_{n,0})w^{2n}+
e_3(w ),
\end{equation}
in which  and also in the rest of discussion, we use $q_j(a_{n,0},\overline
a_{n,0})$ to denote a polynomial which is  symmetric in  $a_{n,0}$ and
$\overline a_{n,0}$,
and also use $e_j(w )$ to denote terms of  order $<2n$ in
$w$, unless it  is otherwise stated.

Write $\varphi^n$ in the form (\ref{eq:3.0}) with $p_n$ in  place of $p$. 
From (\ref{eq:6.1}) and (\ref{eq:7.8}), we see that the linear part of $p_n$ with respect to the
variables $a_n$ and $\overline a_n$ is given by
\begin{equation}\label{eq:7.11}\begin{split}
Lp_n(\xi,\eta)=&i(e^{in\omega}a_{n,0}+\overline
a_{n,0})\xi^n\sum_{k=0}^{n-1}e^{ink\omega}\\
&+e^{-\frac{i}{2}n\omega}(a_{n,0}+\overline a_{n,0})\eta^n\sum_{k=0}^{n-1}e^{-ink\omega}.
\end{split}\end{equation}
From (\ref{eq:4.1}), one gets
\begin{equation*}
\begin{split}Qp_n=
&
\sum_{k=0}^{n-1}\left(Qp\circ(T_1T_2)^k+Lp\circ(T_1T_2)^k
Lp_k\right)\\
& +\sum_{k=0}^{n-1}Lp_k\left(\xi e^{ik\omega}\partial_\xi Lp\circ (T_1T_2)^k-\eta e^{-ik\omega}\partial_\eta Lp\circ(T_1T_2)^k\right).
\end{split}
\end{equation*}
From (\ref{eq:4.1}) and (\ref{eq:7.8}), it follows that for
$\xi=\zeta_j(0)w$ and $\eta=\zeta_j(0)w^{-1}$,
\begin{equation}\label{eq:7.12}
Lp_k(\xi,\eta)=ik(a_{n,0}+\overline
a_{n,0})(\xi^n+\eta^n).
\end{equation}
Hence, we have
$$
Qp_n\bigl(\zeta_j(0)w,\zeta_j(0)w^{-1}\bigr)=n^2a_{n,0}^2\bigl(\zeta_j(0)w\bigr)^{2n}+q_2(a_{n,0},\overline a_{n,0})w^{2n} +e_4(w ).
$$
From (\ref{eq:7.11}), it also follows that
\begin{equation}\label{eq:7.14}
\begin{split}
\partial_\zeta
Lp_n\bigl(\zeta_j(0)w,\zeta_j(0)w^{-1}\bigr)=&-2sn^2\zeta_j^{2s+n-1}(0)a_{n,0}w^n\\ 
 &+c_3(a_{n,0}+ \overline  a_{n,0})w^n+e_5(w ),
\end{split}\end{equation}
in which  $c_3$ is a constant and $e_5(w)$ contains terms of order
$<n$ in $w$.

 By (\ref{eq:4.13}), one gets
$$
\bigl(1+h(\zeta,w)\bigr)^{\frac{1}{2s}}=1+\frac{1}{2sni\zeta^{2s}}p_n(\zeta
w,\zeta w^{-1})+c_4p_n^2(\zeta
w,\zeta w^{-1})+O(p_n^3),
$$ 
in which $c_4$ is a constant. Then for the solution $\zeta=\zeta_j(w, a,\overline a)$, one has
\begin{equation*}\begin{split}
Q(1+h(\zeta,w ))^{\frac{1}{2s}}&=\frac{1}{2sn\zeta_j^{2s}(0)i}\bigl(Qp_n(\xi,
\eta)+\partial_\zeta p_n(\xi,\eta)\cdot L\zeta(\xi,\eta)\bigr)\\ 
&\quad +c_5
L\zeta\cdot Lp_n(\xi,\eta)+
c_4Lp_n^2(\xi,\eta)
\end{split}\end{equation*}
with $\xi=\zeta_j(0)w$ and $ \eta=\zeta_j(0)w^{-1}$. Notice that 
$$
L\zeta=-\zeta_j(0)
\frac{Lp_n\bigl(\zeta_j(0)w,\zeta_j(0)w^{-1}\bigr)}{2sn\zeta_j^{2s}(0)i}=-\zeta_j(0)
\frac{a_{n,0}+\overline a_{n,0}}{2s\zeta_j^{2s}(0)}\zeta_j^n(0)w^n+e_6(w ),
$$
where $e_6(w )$ contains terms of order
$<n$ in $w$. Thus, (\ref{eq:7.14}) gives
$$
\partial_\zeta Lp_n\cdot
L\zeta=n^2a_{n,0}^2\zeta_j^{2n}(0)w^{2n}+q_3(a_{n,0},\overline
a_{n,0})w^{2n}+e_7 (w ).
$$

Now it is easy to see that the solution $\zeta_j(w, a,\overline a)$
to (\ref{eq:4.16}) can be written as
\begin{equation}\label{eq:7.15}
Q\zeta_j
= \frac{in\zeta_j^{2n-2s+1}(0)}{s}\bigl(a_{n,0}^2+h_{n,j}(a,\overline
a)\bigr)w^{2n}+e_{8}(w ),
\end{equation}
in which  $h_{n,j}(a,\overline a)$ is power series converging for
$(a,\overline a)\in D^\infty\times D^\infty$. Furthermore,
  $h_{n,j}(a,\overline a)$ starts with the quadratic terms which is
symmetric in $a_{n,0}$ and $\overline a_{n ,0}$, if $\widetilde a(\xi,\eta)=a_{n,0}\xi^n$.

To complete the proof of Theorem~\ref{thm:7.1}, we put
$$
\zeta_j(w, a,\overline a)=\sum_{k=-\infty}^{\infty}
\zeta_{j,k}(a,\overline a)w^k.
$$
Assume that $\gamma_{n,j}\cap{\Bbb R}^2$ has an accumulation point in
$\Omega$. Then the intersection is a real analytic curve. Notice that $(\xi,\eta)\in {\Bbb R}^2$ impies that  $\zeta$ and $w$ are both
real or pure imaginary.
Hence we get

(a) $\overline\zeta_{j,k}(\overline a,a)=\zeta_{j,k}(a,\overline a),
\quad \mbox{for}\ w=\overline w,\  \mbox{or}$

(b) $\overline\zeta_{j,k}(\overline a,a)=-(-1)^k\zeta_{j,k}(a,\overline a),
\quad \mbox{for}\ w=-\overline w. $

\noindent In particular, we see that 
$\zeta_{j,2n}(a,\overline a)/\zeta_{j,0}(a,\overline a)$ is real.
Since $\zeta_j^{2n-2s}(0)$ is real for $s|n$, then (\ref{eq:7.15}) gives
$$
a_{n,0}^2+\overline a_{n,0}^2+2\widetilde h_{n,j}(a,\overline a)=0,
$$
in which $\widetilde h_{n,j}$ is the real part of $h_{n,j}(a,\overline
a)$. Put
$$
H_n(a,\overline a)=\prod_{j=0}^{2s-1}
\bigl(a_{n,0}^2+\overline a_{n,0}^2+2\widetilde h_{n,j}(a,\overline a)\bigr).
$$
Since $h_{n,j}(a,\overline a)$ is symmetric in $a_{n,0}$ and $\overline
a_{n,0}$
for $a(\xi,\eta)=a_{n,0}\xi^n$, 
then one has $H_n(a,\overline a)\neq 0$. Therefore,
Theorem~\ref{thm:7.1} is proved.

Analogous to Theorem~\ref{cor:1.2}, we have
\begin{cor}\label{cor:7.2}{\it
There is a sequence of power series $H_k(a,\overline a)$ such that
the transformation $\varphi_a,$ defined by 
$(\ref{eq:7.3})$ and $(\ref{eq:7.4}),$
cannot be
transformed into the normal form $(\ref{eq:normalform})$ through any convergent
transformation with  the reality condition$,$  provided 
$H_k(a,\overline a)$ $\neq 0$ for infinitely
many $k.$
}
\end{cor}
For the proof, we keep the notations in Theorem~\ref{thm:7.1}.
Choose a sequence of positive integers
$n_k\ (k=1,2,\ldots)$ such that for each $n=n_k$, the number $\beta$  in
(\ref{eq:beta}) satisfies $-\delta<\beta<0$. We now let $H_k(a,\overline
a)$ be the power series $H_{n_k}(a,\overline
a)$ constructed at the end of the proof of Theorem~\ref{thm:7.1}. From the
normal form (\ref{eq:normalform}), one can see that given  any small deleted neighborhood of the origin in
${\Bbb R}^2$, $\varphi$ has a continuum of periodic points with
period $n_k$ for $k$ large. Therefore,
$\varphi$ cannot be transformed into the normal form
(\ref{eq:normalform}). Finally, without giving details, we mention  that the sequence of
power series $\{H_k(a,\overline
a)\}$ constructed above is indeed functionally independent.

Theorem~\ref{thm:1.3} is a consequence of Corollary~\ref{cor:7.2} and the
intrinsic property of pairs of involutions described early in this section.

\bibliographystyle{plain}

\end{document}